\documentclass{birkjour}
\newtheorem{thm}{Theorem}

\newtheorem{lem}[thm]{Lemma}

\newcommand{\iy}{\infty}

\newcommand{\bC}{{\mathbf C}}
\newcommand{\bT}{{\mathbf T}}
\newcommand{\bZ}{{\mathbf Z}}

\newcommand{\bR}{{\mathbf R}}

\newcommand{\cB}{{\mathcal B}}

\newcommand{\al}{\alpha}
\newcommand{\be}{\beta}
\newcommand{\ga}{\gamma}

\newcommand{\de}{\delta}

\newcommand{\la}{\lambda}
\newcommand{\La}{\Lambda}

\newcommand{\om}{\omega}
\newcommand{\ph}{\varphi}

\newcommand{\vro}{\varrho}
\newcommand{\si}{\sigma}
\newcommand{\tht}{\theta}
\newcommand{\sm}{\setminus}
\newcommand{\n}{\|}

\newcommand{\ti}{\widetilde}
\newcommand{\ov}{\overline}

\newcommand{\ot}{\otimes}

\newcommand{\tia}{{\tilde{a}}}

\begin{document}
	
	%-------------------------------------------------------------------------
	% editorial commands: to be inserted by the editorial office
	%
	%\firstpage{1} \volume{228} \Copyrightyear{2004} \DOI{003-0001}
	%\seriesextraline{This is the Concrete Title of this Book\br H.E. R and S.T.C. W, Eds.}
	%\firstpage{1}

	\title[Spectral radius of the Bell-CHSH operator]
	{The wanted extension\\ of Fujii and Tsurumaru's formula
		for the\\ spectral radius of the Bell-CHSH operator}

	%----------Author 1
	\author[B\"ottcher]{Albrecht B\"ottcher}
	
	\address{%
		Fakult\"at f\"ur Mathematik\\
		TU Chemnitz\\
		09107 Chemnitz\\
		Germany}

	\email{aboettch@mathematik.tu-chemnitz.de}
	
	%\thanks{This work was completed with the support of our
	%----------Author 2
	\author[Spitkovsky]{Ilya M. Spitkovsky}
	\address{Division of Science and Mathematics\\
		New York  University Abu Dhabi (NYUAD)\\
		Saadiyat Island,
		P.O. Box 129188\\
		Abu Dhabi\\ 
		United Arab Emirates}

	\email{ims2@nyu.edu, imspitkovsky@gmail.com}
	
	
	\date{January 12, 2026}
	\dedicatory{In memory of Heinz Langer,\\a prominent operator theorist\\ and a dear friend}
	%%% ----------------------------------------------------------------------
	
	\begin{abstract}
		This paper is motivated by a  recent paper
		of Yuki Fujii and Toyohiro Tsurumaru in which they established a beautiful formula
		for the spectral radius of the Bell-CHSH operator on finite-dimensional Hilbert spaces. 
		To tackle the operator on infinite-dimensional spaces, they elaborated a method based
		on appropriate approximation of commutators of infinite-dimensional orthogonal projections by 
		commutators of orthogonal projections on finite-dimensional spaces. We here give a proof of
		Fujii and Tsurumaru's original formula that works in all dimensions. We also present an alternative
		approximation procedure, uncover the connection of the problem with block Toeplitz operators, 
		and derive good estimates and explicit expressions for the spectral radius in concrete
		cases.
	\end{abstract}
	
	%%% ----------------------------------------------------------------------
	\maketitle
	%%% ----------------------------------------------------------------------
	
	\section{The Bell-CHSH operator}\label{S1}
	
	%\noindent
	Let $P_1,Q_1$ and $P_2,Q_2$ be two pairs of orthogonal projections on separable complex Hilbert spaces
	$V_1$ and $V_2$, respectively, 
	and put $A_i=2P_i -I, B_i=2Q_i-I$. The Bell-CHSH operator is the operator
	\[
	\cB=(A_1+B_1) \otimes A_2+ (A_1-B_1) \otimes B_2
	\]
	on $V_1 \otimes V_2$. The spectral radius $\vro(\cB)$ is of interest in quantum mecanics.
	See, e.g., \cite{CB,  CHSH, FT, GRS, KT, L, T}.
	Taking into account that $A_i$ and $B_i$ are involutions, $A_i^2=B_i^2=I$,
	a simple straightforward
	computation gives
	\[\cB^2=4(I \ot I)+[A_1,B_1] \ot [B_2,A_2],\]
	where $[A,B]=AB-BA$. Consequently, since $\vro(S)=\vro(-S)$ for every operator $S$, we get the
	well-known formula
	\begin{equation}\label{Eq1}
		\vro(\cB)=\sqrt{4+\vro([A_1,B_1]) \vro([A_2,B_2])},
	\end{equation}
	As $[A_i,B_i]=4[P_i,Q_i]$, we are so led to the study of
	the spectral radius of the commutator of two selfadjoint involutions or 
	equivalently, of two orthogonal projections.

	\smallskip
	Thus, let $P,Q$ be two orthogonal projections on a separable complex Hilbert space~$V$
	and put $A=2P-I, B=2Q-I$. We have $\n A \n^2 =\n AA^*\n =\n A^2\n =1$ and, analogously,
	$\n B \n =1$. It follows that 
	\[\vro([A,B]) \le \n [A,B]\n \le \n A \n \,\n B \n +  \n B \n \,\n A \n =2,\]
	and hence (1) implies that $\vro(\cB) \le \sqrt{4+2\cdot 2}=2\sqrt{2}$, which is
	Tsirelson's famous inequality~\cite{T}. We here are interested in expressions for $\vro([P,Q])$
	and thus for $\vro([A,B])=4\vro([P,Q])$ and $\vro(\cB)$ that allow us the exact computation of these
	quantities for concrete operators $P,Q$.
	
	\section{An exact formula for the spectral radius}\label{S2}
	
	For a bounded linear operator $T$, we denote by $\si(T)$ the spectrum of $T$.

	\begin{thm}\label{Theo2}
		We have
		\[\vro([A,B])=\max_{\la \in \si(A+B)}\sqrt{\la^2(4-\la^2)}=\sqrt{\la_0^2(4-\la_0^2)},\]
		where $\la_0$ in $\si(A+B)$ is a point
		for which $\la^2$ is closest to $2$.
	\end{thm}
	
	{\em Proof.}
	We use Halmos' theorem 	\cite{H} as quoted in Theorem 1.2 of~\cite{BS}.
	Let us first assume that $P$ and $Q$ are in general position,
	which means that the terms $(1,1,0,0)$ and $(1,0,1,0)$ in Theorem 1.2 of~\cite{BS}
	are absent.
	Then $A+B-\la I$ is unitarily similar to  $\begin{pmatrix}
		2(I-H)-\la I & 2W\\
		2W & 2(H-I)-\la I
	\end{pmatrix}$. By computing the determinant of this matrix we see that $\la$ is in the spectrum of $A+B$ if and only if
	$1-\la^2/4$ is in the spectrum of $H$. The spectrum of $W^2=H(I-H)$ is therefore
	the set of 
	\[\left(1-\frac{\la^2}{4}\right)\left(1-\left(1-\frac{\la^2}{4}\right)\right)
	=\frac{\la^2}{4}\left(1-\frac{\la^2}{4}\right)\]
	where $\la$ ranges over
	the spectrum of $A+B$. The spectrum of $[P,Q]$ is the set of $\nu$ for which 
	$\det  \begin{pmatrix}
		-\nu I & W\\
		-W & -\nu I
	\end{pmatrix}=\nu^2 I+W^2$
	is not invertible. These values $\nu$ are purely imaginary and we may write $\nu=i\mu$ with real $\mu$.
	Consequently, the spectrum of $[P,Q]$ is the set of $\nu$ such that
	$-\nu^2=\mu^2=\frac{\la^2}{4}\left(1-\frac{\la^2}{4}\right)$ with $\la$ in the spectrum of $A+B$.
	The maximum of $|\mu|$ is attained at a $\la_0$ for which $\la^2/4$ is closest to $1/2$ or equivalently,
	at a $\la_0$ for which $\la^2$ is closest to $2$.
	For this $\la_0$, we have 
	\[\vro([A,B])=4\vro([P,Q])=4|\mu|=4 \sqrt{\frac{\la_0^2}{4}\left(1-\frac{\la_0^2}{4}\right)}
	=\sqrt{\la_0^2(4-\la_0^2)},\]
	as asserted. 
	
	\smallskip
	If $P,Q$ are not in general position, then $[P,Q]$ is unitarily similar
	to 
	\begin{equation}\label{noi}
		(0,0,0,0) \oplus \begin{pmatrix}
			0 & W\\
			-W & 0
		\end{pmatrix}
	\end{equation}
	and hence
	the spectrum of $[P,Q]$
	is the union of $\{0\}$ and the above set of $\nu$ for which 
	$-\nu^2=\frac{\la^2}{4}\left(1-\frac{\la^2}{4}\right)$ with $\la$ in the spectrum of
	$\begin{pmatrix}
		2(I-H) & 2W\\
		2W & 2(H-I)
	\end{pmatrix}$, and we therefore get the same expression for $\vro([A,B])$. \;\: $\square$

	\medskip
	For finite-dimensional spaces $V$, this theorem was established by Fujii and Tsurumaru in~\cite{FT}.
	They wrote
	``when $V$ is infinite-dimensional, it is not known
	whether or how the equality $\vro([A,B])=\sqrt{\la_0^2(4-\la_0^2)}$ can be extended, either in the form of an
	equality or an inequality.'' Theorem~1 is the extension they were asking about.
	
	\smallskip
	The meta theorem of our paper \cite{RS} is that Halmos’ two projections theorem 
	is something like Robert Sheckley’s Answerer: no question
	about two orthogonal projections
	will go unanswered, provided the question is not foolish, and questions that are not foolish
	are those where one already knows most of the answer.
	Fujii and Tsurumaru also employed Halmos' theorem, but they did this only for
	the irreducible $2 \times 2$ blocks of the orthogonal decomposition of $[A,B]$ resulting from
	the unitary representation of the infinite dihedral group $\bZ_2 \ast \bZ_2$ on $V$. We here
	apply Halmos' theorem in its full generality directly to $[A,B]$, without using an orthogonal
	decomposition as an intermediate step. Since we knew the result of~\cite{FT} in advance
	and thus most of the answer, the present study is a nice illustration of the message we tried to 
	convey in~\cite{RS}.

	\smallskip
	Incidentally, and this is again in the vein of~\cite{RS},
	one may also employ Halmos' theorem in its full generality to prove Tsirelson's inequality.
	Indeed, we know that $[P,Q]$ is unitarily similar to the operator (2), it is easily seen that the spectral radius of this operator is equal to the spectral radius of $W$, and the latter
	is $\sqrt{x(1-x)}$ where $x$ is the point in the spectrum of $H$ that is closest
	to $1/2$. Consequently, we always have $\vro([P,Q]) \le 1/2$, whence 
	$\vro([A,B]) \le 2$ and $\vro(\cB) \le 2\sqrt{2}$.  
	
	\smallskip
	In \cite{FT} it is claimed that the $\la_0$ 
	in the formula $\vro([A,B])=\sqrt{\la_0^2(4-\la_0^2)}$
	is the $\la_0$ closest to $\sqrt{2}$. However, there is
	a subtlety in the matter: a $\la_0$ closest to $\sqrt{2}$ need not
	be  a $\la_0$ for which $\la^2$ is closest to $2$. Indeed, if, for example, $A$ and $B$ are given
	on $\bC^3$ by the selfadjoint and involutive matrices
	\[A=\begin{pmatrix}
		1 & 0 & 0\\
		0 & 1 & 0\\
		0 & 0 & -1
	\end{pmatrix}, \quad
	B=\begin{pmatrix}
		1 & 0 & 0\\
		0 & 2x-1 & \sqrt{x(1-x)}\\
		0 & \sqrt{x(1-x)} & 1-2x
	\end{pmatrix}
	\]
	with $x=0.01$, then the spectrum $\si(A+B)$ equals $\{-2\sqrt{x},2\sqrt{x},2\}=\{-0.2,0.2,2\}$, 
	the point $\la=2$ is closest to $\sqrt{2}$, but $\la=\pm 0.2$ are the points for which $\la^2$ is closest
	to $2$.

	\section{Pairs in convergent-angle one-shifted form}\label{S3}
	
	Let $\tht \in (0,\pi)$ and put $s=\sin\tht$ and $c=\cos\tht$. 
	Note that $c^2+s^2=1$ and thus $\frac{1-c}{s}=\frac{s}{1+c}$ and  $\frac{1+c}{s}=\frac{s}{1-c}$.
	Furthermore, if $\tht < \pi/2$ then $c+s >1$, and if $\tht >\pi/2$ then $-1<c+s <1$.
	Fujii and Tsurumaru consider pairs of involutions in so-called {\em constant-angle one-shifted form}. Such a pair is 
	given by
	\begin{eqnarray*}
		& & A(\tht)={\rm diag}\left(\begin{pmatrix} c & s\\ s & -c \end{pmatrix}, \begin{pmatrix} c & s\\ s & -c \end{pmatrix},
		\begin{pmatrix} c & s\\ s & -c \end{pmatrix}, \ldots \right),\\
		& & B(\tht)={\rm diag}\left(1, \begin{pmatrix} c & s\\ s & -c \end{pmatrix}, \begin{pmatrix} c & s\\ s & -c \end{pmatrix},
		\begin{pmatrix} c & s\\ s & -c \end{pmatrix}, \ldots \right).
	\end{eqnarray*}
	In \cite{FT}, it  is claimed that $\si(A(\tht)+B(\tht)) =[-2s,2s]$ for all $\tht \in (0,\pi)$. 
	In fact, the following is true.
	
	\begin{thm}\label{Theo3}
		If $\tht \le \pi/2$ then $\si(A(\tht)+B(\tht)) =[-2s,2s] \cup \{2\}$, and if
		$\tht \ge \pi/2$ then $\si(A(\tht)+B(\tht)) =[-2s,2s]$.
	\end{thm}
	
	We will give two proofs. They are both based on the observation that
	$A(\tht)+B(\tht)=T(a)+K$ where $T(a)$ is the infinite 
	tridiagonal Toeplitz matrix with the symbol $a(e^{i\ph})=s(e^{-i\ph}+e^{i\ph})=2s \cos\ph$ and
	$K$ is the matrix whose $1,1$ entry is $1+c$ and the other entries of which are zero.
	We here encounter an instance of a situation in which work with infinite Toeplitz matrices
	is simpler than working with finite Toeplitz matrices. We note that the appearance of the outlier
	$2$ should not come as too much a surprise. If $\tht=0$, then
	\[A(0)={\rm diag}(1,-1,1,-1,  \ldots),\quad B(0)={\rm diag}(1,1,-1,1,-1,  \ldots),\]
	whence  $A(0)+B(0)={\rm diag}(2,0,0,0,0,\ldots)$ with spectrum $\{0,2\}$. On the other
	hand, if $\tht=\pi$, we get
	\[A(\pi)={\rm diag}(-1,1,-1,1,  \ldots),\quad B(\pi)={\rm diag}(1,-1,1,-1,1,  \ldots)\]
	and thus $A(\pi)+B(\pi)={\rm diag}(0,0,0,0,0,\ldots)$ with spectrum $\{0\}$.
	
	\smallskip
	The essential
	spectrum $\si_{\rm ess}(T)$ is the set of $\la \in \bC$ for which $T-\la I$ is not Fredholm,
	that is, not invertible modulo compact operators.
	
	\medskip
	{\em First proof.} 
	The essential spectrum of $T(a)+K$ is $[-2s,2s]$, and hence the spectrum of $T(a)+K$
	is a subset of $\bR$ containing $[-2s,2s]$. Take $\la \in \bR \sm [-2s,2s]$. Since
	$T(a)+K-\la I$ is Fredholm of index zero, $\la$ is in the spectrum of $T(a)+K$ if and only
	if the kernel of $T(a)+K-\la I$ is not trivial.

	\smallskip
	Put $\mu=\la/s$. The  equation $(T(a)+K-\la I)x=(T(a)+K-s\mu I)x=0$ is the infinite system
	\begin{eqnarray*}
		& & (1+c-s\mu)x_0+sx_1=0\\
		& & x_0-\mu x_1+x_2=0\\
		& & x_1-\mu x_2 +x_3=0\\
		& & \ldots .
	\end{eqnarray*}
	Suppose this system has a non-trivial solution.
	Ignoring the first equation, we get $x_n=\al q_1^n +\be q_2^n$ where $q_1,q_2$ satisfy
	the equation
	$q^2-\mu q +1=0$. As $q_1q_2=1$ and $|q_1+q_2|=|\mu| >2$, we have $|q_1| <1$ and $|q_2| >1$.
	Thus, we get $x_n=\al q_1^n$. Since $q_1-\mu=-1/q_1$, the first equation of the infinite system
	becomes
	\[0=(1+c-s\mu)\al+s\al q_1=\al [1+c+s(q_1-\mu)]=\al[1+c-s/q_1],\]
	and since $\al\neq 0$, it follows that $q_1=s/(1+c)$. 
	
	\smallskip
	We have $s/(1+c)=(1-c)/s \ge 1$  for 
	$\tht \ge \pi/2$, contradicting the requirement that $|q_1| <1$. Consequently, for $\tht \ge \pi/2$ 
	no point outside $[-2s,2s]$ is in $\si(A(\tht)+B(\tht))$.
	
	\smallskip
	So let $\tht < \pi/2$. In that case $q_1=s/(1+c)=(1-c)/s$ lies in $(0,1)$, we get $q_2=1/q_1=(1+c)/s$, 
	and hence
	\[
	\frac{\la}{s}=\mu=q_1+q_2=\frac{1+c}{s}+\frac{s}{1+c}=\frac{2(1+c)}{s(1+c)}=\frac{2}{s}.
	\]
	It results that $\la=2$. For $\la=2$ and $x_n=\al q_1^n$, the first equation of the system reads
	\[(c-1)\al+s\al q_1=\al[c-1+s(1-c)/s]=0,\]
	the remaining equations amount to
	\[1-\frac{2}{s}q_1+q_1^2=1-\frac{2}{s}\frac{s}{1+c}+\frac{1-c}{s}\frac{s}{1+c} =0,\]
	and as this is true, $\la=2$ belongs indeed to the spectrum of $T(a)+K$. \;\: $\square$
	
	\medskip
	{\em Second proof.} We know that $[-2s,2s]$ is contained in the spectrum of $T(a)+K$. 
	Take again $\la \in \bR \sm [-2s,2s]$. The spectrum of $T(a)$ is $[-2s,2s]$ and hence
	$T(a)-\la I=T(a-\la)$ is invertible. So we may write
	\[T(a)+K -\la I=T(a-\la)[I+T^{-1}(a-\la)K],\]
	and it follows that $T(a)+K-\la I$ is invertible if and only if so is $I+T^{-1}(a-\la)K$.
	The matrix of the latter operator is the identity matrix plus $1+c$ times the first column of
	$T^{-1}(a-\la)$, and therefore it is invertible if and only if its $1,1$ entry is nonzero.
	
	\smallskip
	To get the inverse of $T(a-\la)$, we employ Wiener-Hopf factorization. With $t=e^{i\ph}$, we
	have
	\begin{eqnarray*}
		& & a(t)-\la =s(t^{-1}+t)-\la=st^{-1}\left(1-\frac{\la}{s}t+t^2\right)\\
		& & =st^{-1}(t-q_1)(t-q_2)=-sq_2\left(1-\frac{q_1}{t}\right)\left(1-\frac{t}{q_2}\right),
	\end{eqnarray*}
	where $q_1,q_2$ are as in the first proof. It follows that 
	\begin{eqnarray*}
		& & T^{-1}(a-\la)=-\frac{1}{sq_2} T\left((1-t/q_2)^{-1}\right)T\left((1-q_1/t)^{-1}\right)\\
		& & =-\frac{1}{sq_2} \left(\begin{array}{cccc}
			1 & &  & \\
			1/q_2 & 1 & &\\
			1/q_2^2 & 1/q_2 & 1 & \\
			\ldots & \ldots & \ldots & \ldots\end{array}\right)
		\left(\begin{array}{cccc}
			1 & q_1 & q_1^2 & \ldots\\
			& 1 & q_1  &\ldots \\
			&  & 1 & \ldots \\
			\ldots & \ldots & \ldots & \ldots\end{array}\right),
	\end{eqnarray*}
	which shows that the $1,1$ entry of $I+T^{-1}(a-\la)K$ is $1-\frac{1+c}{sq_2}$. This is zero
	if and only if $q_2=(1+c)/s=s/(1-c)$. As $|q_2| >1$ and  $s/(1-c) \le 1$ for $\tht \ge \pi/2$,
	we conclude that $\la$ is not in the spectrum of $T(a)+K$ if $\tht \ge \pi/2$. So let $\tht < \pi/2$.
	Then $q_2 >1$  and since
	$1-\frac{\la}{s}q_2+q_2^2=0$, we obtain
	\[\frac{\la}{s}=\frac{1}{q_2}(1+q_2^2)=\frac{1}{q_2}+q_2=\frac{1+c}{s}+\frac{s}{1+c}=\frac{2(1+c)}{s(1+c)}=\frac{2}{s}.\]
	Thus, necessarily $\la=2$. Conversely, if $\la=2$, then  $q_2=(1+c)/s=s/(1-c)$ satisfies the equation
	$1-\frac{\la}{s}q_2+q_2^2=0$ and is greater than $1$. As  $1-\frac{1+c}{sq_2}=1-\frac{1+c}{s}\frac{s}{1+c}=0$,
	we see that the $1,1$ entry of  $I+T^{-1}(a-\la)K$ is zero, and consequently, $\la=2$ is in the spectrum
	of $T(a)+K$. \;\: $\square$
	
	\medskip
	The following theorem is already in \cite{FT}, but the proof given there is based
	on the incorrect assumption that $\si(A(\tht)+B(\tht)) =[-2s,2s]$ for  $\tht \in (0,\pi/4)$,
	which mistake was fortunately neutralized by the fact that the outlier $\la=2$, though
	being closest to $\sqrt{2}$, is not the $\la$ for which $\la^2$ is closest to $2$.

	\begin{thm}\label{Theo4}
		In the case of a constant-angle one-shifted form we have
		\[
		\vro([A(\tht),B(\tht)])
		= \left\{\begin{array}{lll}
			2\sin(2\tht) &
			\mbox{if} & 0 < \tht \le \pi/4,\\
			2 & \mbox{if} &\pi/4 \le \tht \le 3\pi/4,\\
			2|\sin(2\tht)| &  \mbox{if} & 3\pi/4 \le \tht < \pi. \end{array}\right.
		\]
	\end{thm}
	
	{\em Proof.}
	We use Theorems 1 and 2.
	If  $\tht \le \pi/4$, then $\la^2=2^2=4$ has distance $2$ to $2$ whereas  $\la^2=(2s)^2=4s^2$
	is at the distance $2-4s^2$ to $2$. Consequently, $\la_0=2s$ is a point in the spectrum for which
	$\la^2$ is closest to $2$. This implies that
	\[\vro([A(\tht),B(\tht)])=\sqrt{4s^2(4-4s^2)}=4s\sqrt{1-s^2}=4sc=2\sin(2\tht).\]
	If $\pi/4 \le \tht \le 3\pi/4$,
	then $\sqrt{2} \in [-2s,2s]$ and hence $\la_0=\sqrt{2}$ minimizes the distance between $\la^2$ and $2$.
	Thus, $\vro([A(\tht),B(\tht)])=\sqrt{2(4-2)}=2$. Finally, if $3\pi/4 \le \tht <\pi$,
	then again $\la_0^2=(2s)^2=4s^2$ is closest to $2$, and since $c <0$ in this case, we obtain 
	$\vro([A(\tht),B(\tht)])=2|\sin(2\tht)|$. \;\: $\square$
	
	\medskip
	The general case of two involutions in one-shifted form is given by
	\[
	A={\rm diag}\left(r(\om_1), r(\om_2), r(\om_3),\ldots\right),\quad
	B={\rm diag}\left(1, r(\tht_1),r(\tht_2), r(\tht_3), \ldots \right).
	\]
	where $r(\eta)= \begin{pmatrix} \cos \eta & \sin \eta\\ \sin \eta & -\cos \eta \end{pmatrix}$ and $\{\om_j\}_{=1}^\iy, \{\tht_j\}_{=1}^\iy$
	are sequences of points in $(0,\pi)$.
	We say that two such involutions are in {\em convergent-angle one-shifted form}
	if there is a $\tht \in (0,\pi)$ such that $\om_j \to  \tht$ and $\tht_j\to\tht$ as $j \to \iy$. Let $P_n$ be the projection that sends a sequence
	$\{x_j\}_{j=0}^\iy$ to the sequence $\{x_0, \ldots,x_{n-1},0,0,\ldots\}$. If $T$ is an operator on $\ell^2(\bZ_+)$, then $P_nTP_n$
	may be identified with an $n \times n$ matrix in the natural fashion. The spectrum $\si(P_nTP_n)$ is nothing but the set of the
	eigenvalues of the matrix $P_nTP_n$.
	
	\begin{thm}\label{Theo6}
		If $A$ and $B$ are in convergent-angle one-shifted form, then the spectrum of $P_n(A+B)P_n$
		converges in the Hausdorff metric to the spectrum of $A+B$.
	\end{thm} 
	
	{\em Proof.} In the case at hand, we have $A+B=T(a)+K$ where, as above, $T(a)$ is the infinite tridiagonal Toeplitz matrix
	whose symbol is  $a(e^{i\ph})=2s \cos \ph$ with $s=\sin \tht$ and $K$ is a compact operator (given by a tridiagonal matrix
	whose entries go to zero along each diagonal).
	By Proposition~2.33 and Theorem~~4.16 of~\cite{Uni}, the spectrum of
	$P_n(T(a)+K)P_n$ converges in the Hausdorff metric to $\si(T(a)+K) \cup \si(T(a))$. (Note that $T(a)$ is symmetric,
	which implies that the $T(\ti{a})$ occurring in Theorem~4.16 coincides with $T(a)$.) Since $\si(T(a))=[-2s,2s]$ and
	$\si(T(a)+K)$ contains the essential spectrum $\si_{\rm ess} (T(a))=[-2s,2s]$,  we arrive at the assertion. \;\: $\square$

	\medskip
	\noindent
	{\bf Examples.} Let $\om,\de \in (0,\pi)$ and
	\begin{equation}\label{exx}
		A={\rm diag}\left(r(\om), r(\tht), r(\tht),\ldots\right),\quad
		B={\rm diag}\left(1, r(\tht),r(\tht), r(\tht), \ldots \right).
	\end{equation}
	Put $c=\cos\tht,s=\sin\tht,\ga=\cos\om,\de=\sin\om$. The matrix of $A+B=T(a)+K$ is the tridiagonal
	Toeplitz matrix generated by $a(e^{i\ph})=2s\cos\ph$ with the upper left $2 \times 2$ corner matrix
	$\begin{pmatrix} 0 & s\\ s & 0 \end{pmatrix}$ replaced by  $\begin{pmatrix} 1+\ga & \de\\ \de & c-\ga \end{pmatrix}$.
	We know that $[-2s,2s] \subset \si(A+B)$ and that $A+B-\la I$ is Fredholm of index zero for $\la \notin [-2s,2s]$.
	Hence, a point $\la$ outside $[-2s,2s]$ is in $\si(A+B)$ if and only if $A+B$ has a nontrivial kernel in $\ell^2(\bZ_+)$.
	Since $\n A\n =\n B\n =1$, we also have always $\si(A+B) \subset [-2,2]$. Furthermore, if $\pi/4 \le \tht \le 3\pi/4$,
	then $2s \ge\sqrt{2}$  and hence $\sqrt{2} \in \si(A+B)$, which implies that $\vro([A,B]) =\sqrt{2(4-2)}=2$ in this case.
	
	\smallskip
	The equation $(A+B-\la I)x=0$ is the system
	\begin{eqnarray*}
		& & (1+\ga-\la)x_0+\de x_1=0,\\
		& & \de x_0+(c-\ga-\la) x_1+ sx_2=0,\\
		& & sx_n -\la x_{n+1} +sx_{n+2}=0 \;\: (n \ge 1).
	\end{eqnarray*} 
	This is satisfied if and only if $x_n=\al q^{n-1}$ for $n \ge 1$, where $\al \neq 0$ and $q$ is the solution of the equation
	$q^2-(\la/s)q+1=0$ whose modulus is less than $1$, and if $x_0$ can be chosen so that the first two equations
	of the system hold. Since $x_1=1$ and $x_2=q$, the latter is possible if and only if 
	$\de/(1+\ga-\la)=(c-\ga-\la+sq)/\de$, that is,
	\begin{equation} \label{la2}
		\la^2-(c+sq+1)\la +(1+\ga)(c+sq-1)=0.
	\end{equation} 
	
	%%%%%%%%%%%%%%%%%%%%%%%%%%%%%%%%%%%%%%%%%%%%
	%%%%%%%%%%%%%%%%%%%%%%%%%%%%%%%%%%%%%%%%%%%%
	\begin{figure}[tb]
		\begin{center}
			{\includegraphics[height=4cm]{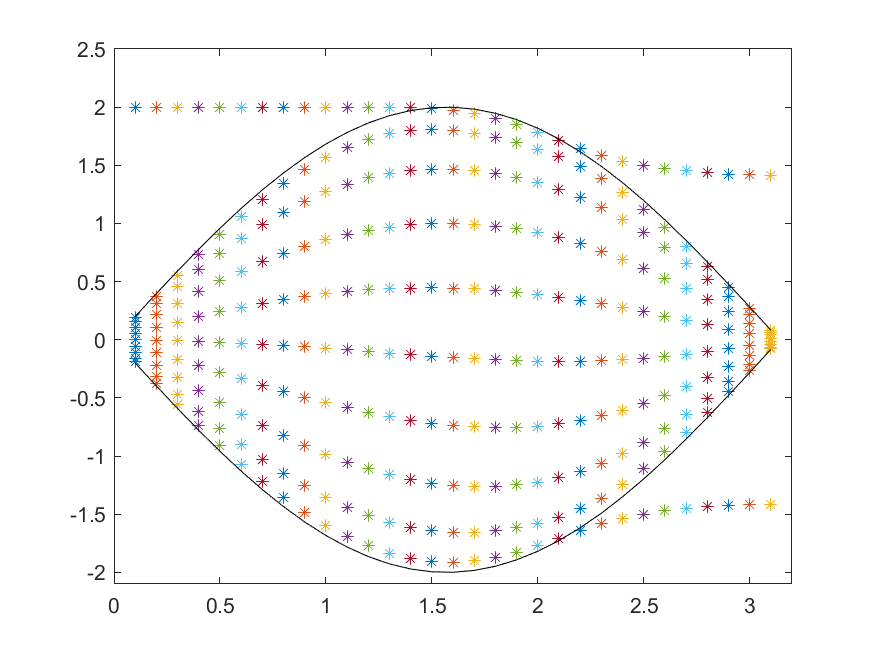}}
			{\includegraphics[height=4cm]{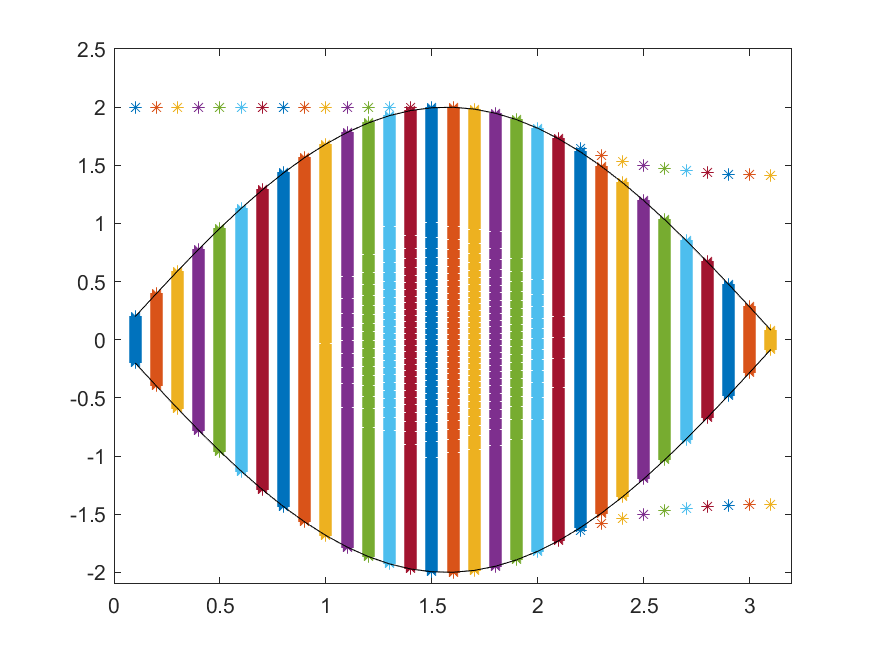}}
			\caption{The spectra of $P_n(A+B)P_n$ multiplied by $i$ for $n=10$ (left) and $n=600$ (right)
				with $A,B$ as in~(3)
				and with $\om=\pi/2$ and 
				$\tht=0.1:0.1:3.1$. The solid curves are plots of $\pm 2\sin\tht$.} 
			\label{Fig1}
		\end{center}
	\end{figure}
	%%%%%%%%%%%%%%%%%%%%%%%%%%%%%%%%%%%%%%%%%%%%
	%%%%%%%%%%%%%%%%%%%%%%%%%%%%%%%%%%%%%%%%%%%%
	
	%%%%%%%%%%%%%%%%%%%%%%%%%%%%%%%%%%%%%%%%%%%%
	%%%%%%%%%%%%%%%%%%%%%%%%%%%%%%%%%%%%%%%%%%%%
	\begin{figure}[tb]
		\begin{center}
			{\includegraphics[height=4cm]{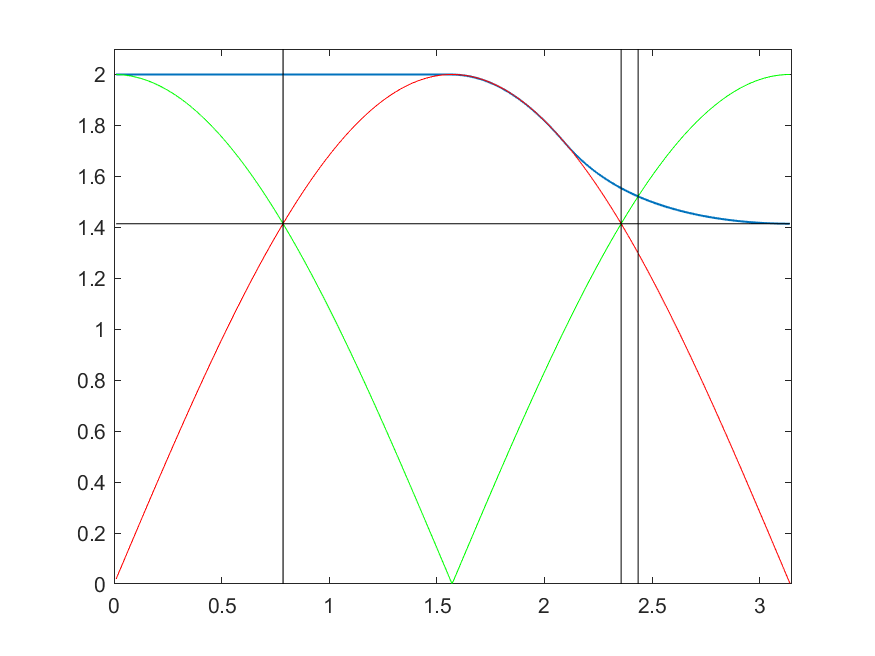}}
			{\includegraphics[height=4cm]{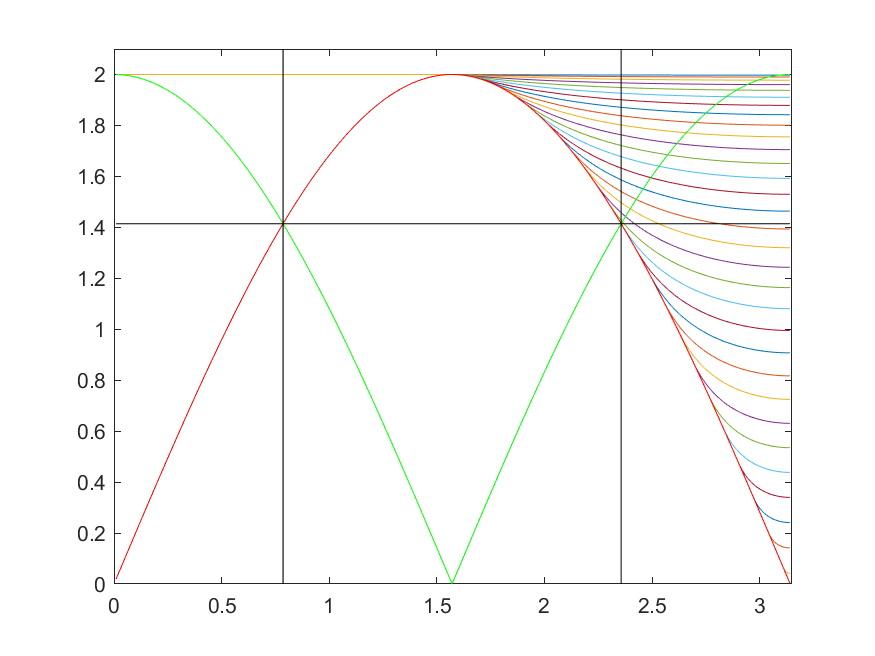}}
			\caption{Plots of $\la_{\max}(\tht)$ where $\la_{\max}(\tht)$ is the largest eigenvalue
				of $P_n(A+B)P_n$ for $n=100$ with $A,B$ as in~(3). In the left picture we took $\om=\pi/2$, 
				in the right picture
				we have $\om=0.1:0.1:3.1$. The red curve is $2\sin \tht$, the green curve is
				$2|\cos \tht|$. The horizontal lines are at the height $\sqrt{2}$. The rightmost
				vertical line in the left picture has the abscissa $x_3=2.4352$.} 
			\label{Fig2}
		\end{center}
	\end{figure}
	%%%%%%%%%%%%%%%%%%%%%%%%%%%%%%%%%%%%%%%%%%%%
	%%%%%%%%%%%%%%%%%%%%%%%%%%%%%%%%%%%%%%%%%%%%
	
	Instead of solving this equation, let us use Theorem~4. Figure~1 shows the result for $\om=\pi/2$
	and $\tht=0.1:0.1:3.1$. To be sure that the outliers are not produced by numerical instabilities, let us check 
	whether~(4) is true for the rational numbers $c=-0.8$ and $s=0.6$, which corresponds to $\tht=\arccos(-0.8) \approx 2.4981$
	and where exact computations can be done by hand.
	The picture suggests that in this case the two outliers are $\la=\pm 1.5$. The roots of $q^2-(1.5/0.6)q+1$ are
	$2$ and $1/2$. Hence, we have to take $q=0.5$. We have $\ga=0$ and $\de=1$ for $\om=\pi/2$.
	Inserting these values in the left-hand side of~(4), we obtain indeed $0$. The same happens for $\la=-1.5$.
	Thus, we may conclude that $\si(A+B)=[-1.2,1.2]\cup\{-1.5,1.5\}$ for $\om=\pi/2$ and $\tht \approx 2.4981$.
	As $1.5^2-2 = 0.25 < 0.56=2-4\cdot 0.6^2$, we see that in the case at hand the spectral radius $\vro([A,B])$ is
	not determined by $s$ but rather by the outliers $\la=\pm 1.5$. Theorem~1 gives 
	$\vro([A,B])=\sqrt{1.5^2(4-1.5^2)} \approx 1.9843$. The value of $\sqrt{4s^2(4-4s^2)}$ is $1.9200$.

	%%%%%%%%%%%%%%%%%%%%%%%%%%%%%%%%%%%%%%%%%%%%
	%%%%%%%%%%%%%%%%%%%%%%%%%%%%%%%%%%%%%%%%%%%%
	\begin{figure}[tb]
		\begin{center}
			{\includegraphics[height=6cm]{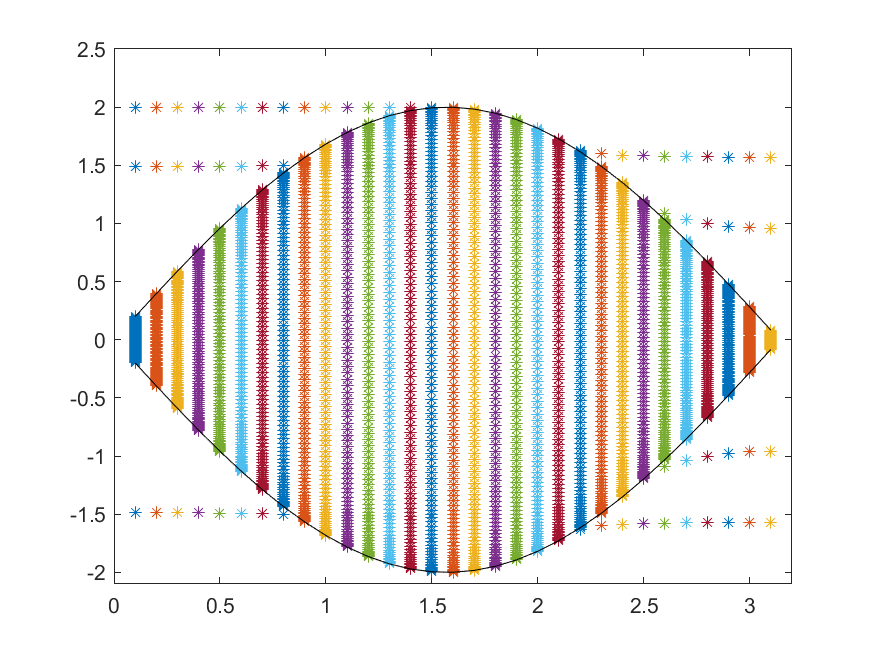}}
			\caption{The spectra of $P_n(A+B)P_n$ multiplied by $i$ for $n=100$
				with $A,B$ as in~(5) and 
				$\tht=0.1:0.1:3.1$.} 
			\label{Fig3}
		\end{center}
	\end{figure}
	%%%%%%%%%%%%%%%%%%%%%%%%%%%%%%%%%%%%%%%%%%%%
	%%%%%%%%%%%%%%%%%%%%%%%%%%%%%%%%%%%%%%%%%%%%

	\smallskip
	The blue curve in the left picture of Figure~2 shows the outliers $\la >0$ for $\om=\pi/2$  
	and $0 <\tht < \pi$. The curve is of constant height $2$ for 
	$0<\tht \le \pi/2$ and then goes down for $\pi/2 < \tht <\pi$. The abscissas of the intersections of
	the red curve $\tht \mapsto 2\sin\tht$ and the green curve $\tht \mapsto |2\cos\tht|$ are $x_1=\pi/4$ and 
	$x_2=3\pi/4 \approx 2.3562$,
	and the abscissa of the intersection of
	the green curve and the blue outlier curve $\tht \mapsto  \la_{\max}(\tht)$ is $x_3 \approx 2.4352$. For $\tht >x_2$,
	we have $\la > \sqrt{2}$ and $2s < \sqrt{2}$,  i.e., both $\la^2-2$ and $2-4s^2$ are positive. 
	It follows that $\la^2-2 < 2-4s^2$ if and only if $\la^2 < 4(1-s^2)=4c^2$. The latter is the case
	exactly if the blue curve is below the green curve, that is, for $\tht > x_3$. In summary, 
	taking for each $\tht$ the
	value of $\la_0 \in \si(A+B)$ for which $\la^2$ is closest to $2$, we arrive at the following: 
	\[\vro([A,B])=\left\{
	\begin{array}{lll} 
		\sqrt{4s^2(4-4s^2)}=2\sin(2\tht) & \mbox{for} & 0 < \tht \le x_1,\\
		\sqrt{2(4-2)}=2 & \mbox{for} & x_1 \le \tht \le x_2,\\
		\sqrt{4s^2(4-4s^2)}=2|\sin(2\tht)| & \mbox{for} & x_2 \le \tht \le x_3,\\
		\sqrt{\la_{\max}(\tht)^2(4-\la_{\max}(\tht)^2)}& \mbox{for} & x_3 \le \tht < \pi.
	\end{array}
	\right.\]
	The right picture of Figure~2 shows the situation for $\om=0.1:0.1:3.1$. Figure~3 is for 
	\begin{equation}\label{pert}
		A={\rm diag}(r(3/2), r(2), r(\tht), r(\tht), \ldots), B={\rm diag}(1,r(5/2), r(\tht), r(\tht),\ldots).
	\end{equation}

	\section{Pairs in one-shifted form with two converging angles} \label{S4}
	
	A natural generalization of pairs in convergent-angle one-shifted form is pairs
	of the form
	\[
	A={\rm diag}\left(r(\om_1), r(\om_2), r(\om_3),\ldots\right),\quad
	B={\rm diag}\left(1, r(\tht_1),r(\tht_2), r(\tht_3), \ldots \right)
	\]
	with numbers $\om,\tht \in (0,\pi)$ such that $\om_j \to \om$ and $\tht_j \to \tht$ as $j \to \iy$.
	Then $A+B=T(a)+K$ with a selfadjoint block Toeplitz matrix $T(a)$ and a compact operator $K$.
	The diagonal block of $T(a)$ is 
	$\begin{pmatrix} \ga-c & \de\\ \de & c-\ga\end{pmatrix}$,
	the block in the subdiagonal equals $\begin{pmatrix} 0 & s\\ 0 & 0\end{pmatrix}$,
	the block in the superdiagonal is equal to
	$\begin{pmatrix} 0 & 0\\ s & 0\end{pmatrix}$, and the other blocks are zero.
	Consequently,
	\begin{equation} \label{sym}
		a(e^{i\ph})=\begin{pmatrix}
			\ga-c &  \de +s e^{i\ph}\\
			\de +s e^{-i\ph} & c-\ga
		\end{pmatrix},
	\end{equation}
	where $\ga=\cos \om, \de=\sin\om, c=\cos \tht, s=\sin \tht$. Put
	\begin{eqnarray*}
		& & \la_1=\sqrt{(c-\ga)^2+(s-\de)^2}=2\left|\sin\frac{\tht-\om}{2}\right|,\\
		& & \la_2=\sqrt{(c-\ga)^2+(s+\de)^2}=2\left|\sin\frac{\tht+\om}{2}\right|.
	\end{eqnarray*}
	Note that $s$ and $\de$ are positive.
	Unless $s=\de$, which case was disposed of in the previous section, we have
	\begin{equation}\label{gapp}
		-2  \le -\la_2 < -\la_1 < c-\ga < \la_1 < \la_2 \le 2.
	\end{equation}
	
	\smallskip
	In cases in which the spectrum of $A+B=T(a)+K$ is available, Theorem~1 gives the spectral
	radius. To determine this spectrum, it would be desirable to have a result like Theorem~4.
	Since the blocks of $T(a)$ are of size $2$, we consider
	only $2n \times 2n$ truncations.   
	It is known that
	\begin{equation}\label{Lalim}
		\si(P_{2n}(T(a)+K)P_{2n}) \to \si(T(a)+K) \cup \si(T(\tia))=: \La
	\end{equation}
	in the Hausdorff metric; see page 197 of \cite{BS}. Here $\tia$ is given by
	\begin{equation} \label{tisym}
		\tia(e^{i\ph})=a(e^{-i\ph})=\begin{pmatrix}
			\ga-c &  \de +s e^{-i\ph}\\
			\de +s e^{i\ph} & c-\ga
		\end{pmatrix}.
	\end{equation}
	Thus, the convergence of $\si(P_{2n}(T(a)+K)P_{2n})$ to some limit set $\La$ is guaranteed.
	However, the appearance of the spectrum $\si(T(\tia))$ in (8) causes some trouble
	since it forces us to find out which parts of the limit set $\La$ constitute $\si(T(a)+K)=\si(A+B)$
	and which parts result only from the contribution of  $\si(T(\tia))$.
	
	\begin{lem} \label{Lem}
		Let $a$ be given by {\rm (6)}. The essential spectra of $T(a)$ and $T(\tia)$ are
		\[\si_{\rm ess}(T(a))=\si_{\rm ess}(T(\tia))=[-\la_2,-\la_1] \cup [\la_1,\la_2]\]
		and the spectrum of $T(\tia)$ is $\si_{\rm ess}(T(\tia))$ if $s \le \de$ and 
		$\si_{\rm ess}(T(\tia)) \cup \{c-\ga\}$ if $s > \de$.
	\end{lem}
	
	{\em Proof.} The essential spectrum of $T(a)$ is the set of $\la \in \bR$ for which the determinant
	of $a(e^{i\ph})-\la I$ vanishes for some $e^{i\ph} \in \bT$; see, e.g., Theorem~6.5 of~\cite{BS}.
	We have
	\[\det (a(e^{i\ph})-\la I)=\la^2 -(c-\ga)^2 -\de^2 -2\de s \cos \ph -s^2=0\]
	if and only if
	$\la^2=(c-\ga)^2 +\de^2 +2\de s \cos \ph +s^2$, and since $2\cos \ph$ fills the segment $[-2,2]$,
	it follows that the essential spectrum of $T(a)$ is equal to $[-\la_2,-\la_1] \cup [\la_1,\la_2]$.
	As $\cos\ph=\cos(-\ph)$, we get the same result for $\si_{\rm ess}(T(\tia))$.
	
	\smallskip
	That the spectrum $\si(T(\tia))$ is a subset of the segment $[-\la_2, \la_2]$
	can be deduced from general results on
	positively definite Hermitian matrix function; see Section 7.3
	of~\cite{LS}. However, we have also to look whether or not the spectrum contains points
	in the gap $(-\la_1,\la_2)$. We therefore proceed as follows.
	
	\smallskip
	Take $\la \notin \si_{\rm ess}(T(\tia))$. Then the selfadjoint operator $T(\tia)-\la I$ is Fredholm
	of index zero, and hence  $\la$ belongs to $\si(T(\tia))$ if and only if the null space of this
	operator
	is nontrivial. The equation \[(T(\tia)-\la I)\begin{pmatrix} f_+\\ g_+\end{pmatrix}=0\] with $f_+,g_+$ in
	the Hardy space $H^2(\bT)$ is equivalent to the existence of functions $h_-,k_- \in \overline{H^2}(\bT)$
	such that $h_-(0)=k_-(0)=0$ and
	\begin{eqnarray*}
		& & (\ga-c-\la)f_+(z)+(\de+s\ov{z})g_+(z)=h_-(z),\\
		& & (\de+s{z})f_+(z)+(c-\ga-\la)g_+(z)=k_-(z)
	\end{eqnarray*}
	for $z \in \bT$. A moment's thought reveals that this can only be true for $k_-(z)=0$ and $h_-(z)=\al \ov{z}$
	with some constant $\al$. Thus, we arrive at the system
	\begin{eqnarray*}
		& & (\ga-c-\la)zf_+(z)+(\de z+s)g_+(z)=\al,\\
		& & (\de +s z)f_+(z)+(c-\ga-\la)g_+(z)=0.
	\end{eqnarray*}
	Let $\la \neq c-\ga$. Solving the second of these equation for $g_+(z)$ and inserting the result
	in the first equation, we get after some straightforward computation the equation
	$p_+(z)f_+(z)/(c-\ga+\la)=\al$ with
	\[p_+(z)=s\de\left[1+\frac{2-2c\ga-\la^2}{s\de}z+z^2\right].\]
	This equation has a nontrivial solution $f_+ \in H^2(\bT)$ if and only if $1/p_+$ is in $H^2(\bT)$, which in turn happens
	if and only if both roots of $p_+(z)$ lie outside the closed unit disk. As, by Vieta, the product of the two roots
	is $1$, this is impossible. Consequently, a point $\la \neq c-\ga$ never belongs to the spectrum of $T(\tia)$.
	
	\smallskip
	So consider $\la=c-\ga$. Then the system is reduced to 
	\begin{eqnarray*}
		& & 2(\ga-c-\la)zf_+(z)+(\de z+s)g_+(z)=\al,\\
		& & (\de +s z)f_+(z)=0,
	\end{eqnarray*}
	and this holds with a
	nontrivial function $g_+ \in H^2(\bT)$ if and only if $\de z+s$ has no zero in the closed unit disk, that is, if and
	only if $s > \de$. \;\: $\square$

	\medskip
	As said, the point $\la=c-\ga$ causes some trouble.
	We have at least
	the following.
	
	\begin{thm} \label{Theo8}
		Let  $\la_0$ be a point in $\La$ for which $\la^2$ is closest to
		$2$. If $\la_0 \neq c-\ga$, then $\vro([A,B])=\sqrt{\la_0^2(4-\la_0^2)}$, and if 
		$\la_0 = c-\ga$, then 
		\[ 
		\sqrt{\la_*^2(4-\la_*^2)}\le \vro([A,B]) \le \sqrt{\la_0^2(4-\la_0^2)}
		\]
		where $\la_*$ is a point in $\La \setminus \{c-\ga\}$ for which $\la^2$ is closest to $2$.
	\end{thm}
	
	{\em Proof.} 
	The essential spectrum of $T(a)+K$ is equal to the essential spectrum of $T(a)$
	and is a subset of $\si(T(a)+K)$. From Lemma~5 we therefore see that
	\[\si(T(a)+K)=[-\la_2,-\la_1] \cup [\la_1,\la_2] \cup S,\]
	with some (possibly empty) set $S$. Hence, again by Lemma~5,  
	\[\La= [-\la_2,-\la_1] \cup [\la_1,\la_2] \cup S \cup X,\]
	where $X=\{c-\ga\}$ if $c-\ga \in \si(T(a))$ and $X=\emptyset$ if
	$c-\ga \notin \si(T(a))$. Thus, if $\la_0 \neq c-\ga$, then $\la_0 \in \si(T(a)+K)$
	and Theorem~1 gives $\vro([A,B])=\sqrt{\la_0^2(4-\la_0^2)}$, while 
	if $\la_0 =c-\ga$, it may be that $\la_0$ is not in $\si(T(a)+K)$, so that we must
	exclude $\la_0$ from the candidates and have to restrict ourselves to $\La$
	minus $\{\la_0\}$, which gives a lower bound for $\vro([A,B])$. \;\: $\square$

	\medskip
	The following theorem contains two constellations in which $c-\ga$ is ruled out automatically.
	We remark that 
	\begin{eqnarray*}
		& & \sqrt{2} \in (\la_2,2) \;\: \Longleftrightarrow \;\:\tht+\om \in (0,\pi/2) \cup (3\pi/2,2\pi),\\
		& & \sqrt{2} \in (0,\la_1) \;\: \Longleftrightarrow \;\:\tht-\om \in (-\pi/2 ,\pi/2).
	\end{eqnarray*}

	\begin{thm} \label{Theo7}
		Let $\La$ be the limit of $\si(P_{2n}(A+B)P_{2n})$ as $n \to \iy$.

		{\rm (a)} If $\sqrt{2} \in (\la_2,2)$, then 
		$\vro([A,B])=\sqrt{\la_0^2(4-\la_0^2)}$
		where $\la_0 \in \La$
		is a point for which $\la^2$ is closest to $2$. 
		
		{\rm (b)} If $\sqrt{2} \in [\la_1,\la_2]$, then $\vro([A,B])=2$. 
		
		{\rm (c)}  If $\sqrt{2} \in (0, \la_1)$,
		then  $\sqrt{\la_1^2(4-\la_1^2)} \le  \vro([A,B]) \le 2$.
	\end{thm}
	
	{\em Proof.}
	(a) Lemma~5 and (7) tell us that
	$\si(T(\tia))$ is a subset of the segment $[-\la_2,\la_2]$
	and
	that
	the set $\La$ differs from $\si(T(a)+K)$ only by the possible point $c-\ga$ in $\si(T(\tia))$,
	which lies in $[-\la_1,\la_1]$. If $c-\ga=\pm \la_1$ (which happens only for $s=\de$ and thus $\la_1=0$), 
	then $c-\ga$ belongs to $\si(T(a)+K)$.
	Thus, let $c-\ga \in (-\la_1,\la_1)$. Since $\la_2 < \sqrt{2}$,
	we have \[2-\la^2 > 2-\la_1^2 \ge 2-\la_2^2 >0\] for $\la \in (-\la_1,\la_1)$, and hence the point
	$\la=c-\ga$ cannot be a point in $\La$ for which $\la^2$ is closest to $2$.
	It follows that the points $\la$ in $\La$ for which $\la^2$ is closest to $2$
	belong to $\si(T(a)+K)$, and so Theorem~1 gives the assertion.
	
	\smallskip
	(b) If $\sqrt{2} \in [\la_1,\la_2]$, then 
	$\la=\sqrt{2}$ belongs to $\si(T(a)+K)$, it is a point for which $\la^2$ is closest to $2$,
	and from Theorem~1 we get $\vro([A,B])=2$.
	
	\smallskip
	(c) In this case it may be that $\la=c-\ga$ is a point for which $\la^2$ is closest to $2$
	but that this point lies in $\La \setminus \si(T(a)+K)$. 
	Since $\la_1$ is definitely a point in $\si(T(a)+K)$,  
	Theorem~1 implies at least the lower bound $\sqrt{\la_1^2(4-\la_1^2)}$ for $\vro([A,B])$.
	Finally, we know that always $\vro([A,B]) \le 2$. \;\: $\square$
	
	\section{Pairs in one-shifted form with two constant angles} \label{S5}
	
	Let us finally consider the case were both angles are constant,
	\begin{equation} \label{con}
		A={\rm diag}\left(r(\om), r(\om), r(\om),\ldots\right),\;\:
		B={\rm diag}\left(1, r(\tht),r(\tht), r(\tht), \ldots \right)
	\end{equation}
	with two numbers $\om,\tht$ in the interval $(0,\pi)$. As above, let
	\[c=\cos\tht,\; s=\sin\tht,\; \ga=\cos\om,\; \de=\sin\om.\]
	In that case $A+B=T(a)+K$ with $a$ as in~(6) and with
	\[K={\rm diag}(1+c,0,0,\ldots).\] Figure~4 shows the limit sets for $\om=0.3$ and $\om=2.4$.
	These figures suggest that $c-\ga$ belongs to $\si(A+B)$ if and only if $\tht=\om$ or
	$\tht=\pi-\om$, which happens if and only if $c=\ga$ or $c=-\ga$. If $\tht=\om$, then indeed
	$c-\ga=0= \la_1$ and if $\tht=\pi-\om$, then $c-\ga=2c=\pm\la_1$, and $\pm \la_1$ belong to $\si(T(a))+K)$.
	The case $c \neq \pm \ga$ is settled by the following lemma.

	\begin{lem} \label{Lem2} If $A,B$ are given by {\rm (10)} and $\tht\notin\{\om,\pi-\om\}$, then 
		$c-\ga \notin \si(A+B)$.
	\end{lem}
	
	{\em Proof.}
	Assume the contrary: $c-\ga \in \si(A+B)=\si(T(a)+K)$. As $c\neq \ga$, we have $\la_1 >0$.
	Since $T(a)+K-\la I$ is Fredholm
	of index zero for $\la \in (-\la_1, \la_1)$, the equation
	\[(T(a)+K-\la I)
	\begin{pmatrix} f_+\\ g_+\end{pmatrix}=0\] has a nontrivial solution with $f_+,g_+$ in
	$H^2(\bT)$ if $\la=c-\ga$. As in the proof of Lemma~5, this equation can be shown to be
	equivalent to the existence of a constant $\al$ such that
	\begin{eqnarray*}
		& & (1+c)[f_+]_0 +(\ga-c-\la)f_+(z)+(\de+sz)g_+(z)=0,\\
		& & (\de z+s)f_+(z)+(c-\ga-\la)zg_+(z)=\al,
	\end{eqnarray*}
	where $[f_+]_0$ denotes te $0$th Fourier coefficient of $f_+$. For $\la=c-\ga$, this becomes
	\begin{eqnarray*}
		& & (1+c)[f_+]_0 +2(\ga-c)f_+(z)+(\de+sz)g_+(z)=0,\\
		& & (\de z+s)f_+(z)=\al.
	\end{eqnarray*}
	If $s \le \de$, then $\de z+s$ has a zero in the closed unit disk and therefore the second equation
	of the system
	implies that $\al=0$ and that $f_+(z)=0$ for all $z$. Inserting this in the first equation, we see
	that $g_+$ must also vanish identically. But this contradicts our assumption, according to which the
	system must have a nontrivial solution $(f_+,g_+)$.
	
	\smallskip 
	We are left with $s >\de$. The solution $f_+(z)=\al/(s+\de z)$ of the second of these equations is in $H^2(\bT)$.
	We have $[f_+]_0=\al/s$, and hence the first equation takes the form
	\[\frac{(1+c)\al}{s}+\frac{2(\ga-c)\al}{s+\de z}+(\de+s z)g_+(z)=0.\]
	It follows that
	\[g_+(z)=-\frac{\al}{\de+s z}\left[\frac{1+c}{s}+\frac{2(\ga-c)}{s+\de z}\right].\]
	The denominator $\de+s z$ has the zero $z=-\de/s$ in the unit disk. Consequently, for $g_+$ to be
	in $H^2(\bT)$ it is necessary that the term in brackets also vanishes at $z=-\de/s$. Thus, we
	must have
	\[0=\frac{1+c}{s}+\frac{2(\ga-c)}{s+\de (-\de/s)}=\frac{1+c}{s}+\frac{2(\ga-c)s}{s^2-\de^2},\]
	or equivalently, $(1+c)(s^2-\de^2)=2s^2(c-\ga)$. Since $s^2=1-c^2$ and $\de^2=1-\ga^2$, we obtain
	\[(1+c)(\ga^2-c^2)=2(1-c^2)(c-\ga),\]
	and as $c\neq -\ga$, we get $-(\ga+c)=2(1-c)$ and thus $c-\ga=2$. 
	This contradicts (7) and so proves our claim. \;\: $\square$
	
	\smallskip
	Recall that, by Lemma~5,
	\begin{eqnarray*}
		& & s \le \de \;\:\Longleftrightarrow\;\:\tht \in (0,\om] \cup [\pi-\om,\pi)
		\;\:\Longleftrightarrow\;\: c-\ga \notin \si(T(\tia)),\\
		& &  s > \de \;\:\Longleftrightarrow\;\:\tht \in (\om,\pi-\om)
		\;\:\Longleftrightarrow\;\: c-\ga \in \si(T(\tia)),
	\end{eqnarray*}
	
	\begin{thm} \label{Theo9} If $A,B$ are given by {\rm (10)}, then $\si(A+B)=\La$
		for $s \le \de$ and $\si(A+B)=\La \setminus \{c-\ga\}$ for $s >\de$.
	\end{thm}
	
	{\em Proof.} We have $\La=\si(A+B) \cup \si(T(\tia))$ and $A+B=T(a)+K$.
	If $s \le \de$, then Lemma~5 shows that $\si(A+B)=\La$. For $s >\de$,
	we have, again by Lemma~5,  $\si(A+B)=\La$ if $c-\ga  \in \si(T(a)+K)$
	and   $\si(A+B)=\La \setminus \{c-\ga\}$ if $c-\ga  \notin \si(T(a)+K)$.
	But Lemma~8 says that $c-\ga$ is never in $\si(T(a))+K)$ for $s >\de$,
	which implies that  $\si(A+B)=\La \setminus \{c-\ga\}$ in this case. \;\: $\square$

	\begin{thm} \label{Theo10} Let $A,B$ be given by {\rm (10)}. If $0<\om \le \pi/2$, then 
		\[\vro([A,B])=\left\{ \begin{array}{lll}
			2\sin(\tht+\om) & \mbox{for} & 0< \tht \le \pi/2-\om,\\
			2 & \mbox{for} & \pi/2-\om \le \tht \le \pi/2+\om,\\
			2|\sin(\tht-\om)| & \mbox{for} & \pi/2+\om \le \tht < \pi,
		\end{array}\right.\]
		and if $\pi/2 \le \om < \pi$, then
		\[\vro([A,B])=\left\{ \begin{array}{lll}
			2|\sin(\tht-\om)| & \mbox{for} & 0< \tht \le \om-\pi/2,\\
			2 & \mbox{for} & \om-\pi/2 \le \tht \le 3\pi/2-\om,\\
			2\sin(\tht+\om) & \mbox{for} & 3\pi/2-\om \le \tht < \pi.
		\end{array}\right.\]
	\end{thm}
	
	%%%%%%%%%%%%%%%%%%%%%%%%%%%%%%%%%%%%%%%%%%%%
	%%%%%%%%%%%%%%%%%%%%%%%%%%%%%%%%%%%%%%%%%%%%
	\begin{figure}[tb]
		\begin{center}
			{\includegraphics[height=4cm]{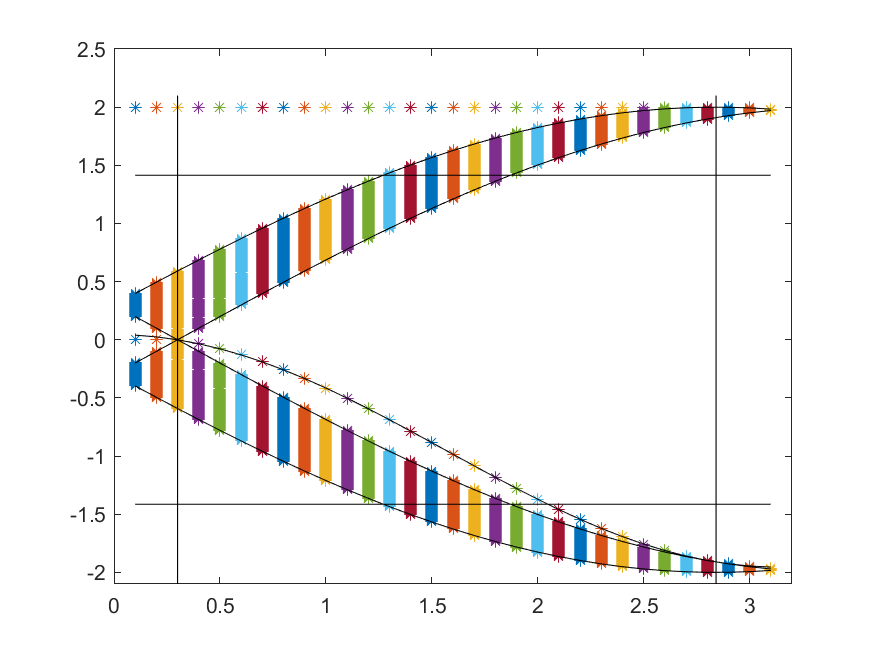}}
			{\includegraphics[height=4cm]{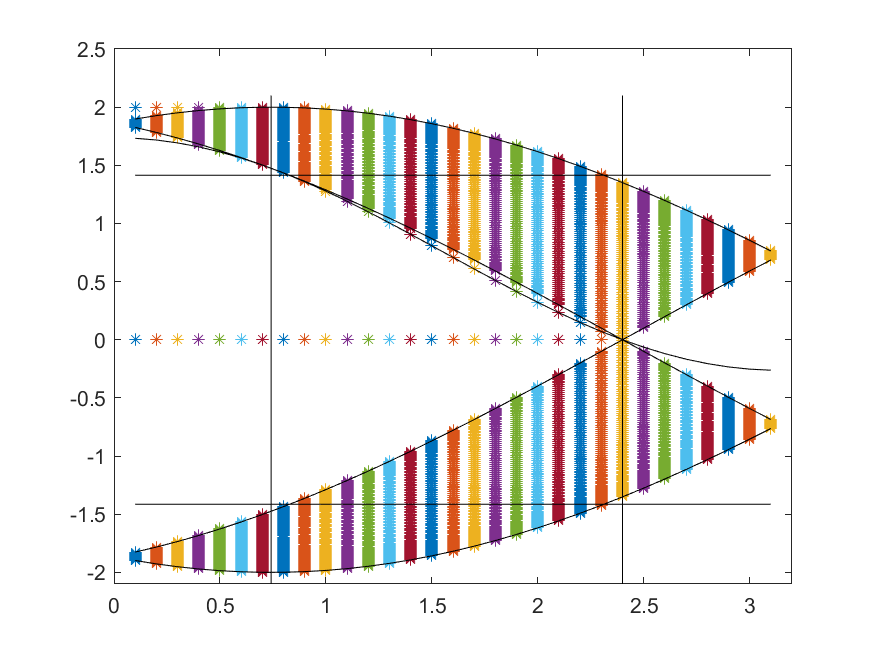}}
			\caption{The spectra of $P_{200}(A+B)P_{200}$ multiplied by~$i$ with 
				$A,B$ as in~(10)
				for $\om=0.3$ (left), $\om=2.4$ (right), and 
				$\tht=0.1:0.1:3.1$. The solid lines plot $c-\ga$, the
				boundaries of $[-\la_2,-\la_1]$ and $[\la_1,\la_2]$,
				and the values $\pm\sqrt{2}$. The vertical lines have the abscissas
				$\om$ and $\pi-\om$.} 
			\label{Fig4}
		\end{center}
	\end{figure}
	%%%%%%%%%%%%%%%%%%%%%%%%%%%%%%%%%%%%%%%%%%%%
	%%%%%%%%%%%%%%%%%%%%%%%%%%%%%%%%%%%%%%%%%%%%
	
	{\em Proof.} Suppose $0 < \om < \pi/2$. In this case the situation is as in the left
	picture of Figure~4. The curve described by $\la_2$ intersects the horizontal line at
	height $\sqrt{2}$ for the $\tht$ given by $2\sin((\tht+\om)/2)=\sqrt{2}$, that is,
	for $\tht=\pi/2-\om$. It is obvious that for $0< \tht <\pi/2-\om$ the point $\la_0=\la_2$
	is a point in $\La$ whose square is closest to $2$, and we know that $\la_0=\la_2 \in \si(A+B)$.
	Thus, for $0 < \tht < \pi/2-\om$ we get
	\[\vro([A,B])=\sqrt{\la_2^2(4\!-\!\la_2^2)}=\sqrt{4\sin^2\!\frac{\tht\!+\!\om}{2}\!\left(\!4\!-\!4\sin^2\!\frac{\tht\!+\!\om}{2}\right)}
	=2\sin(\tht+\om).\]
	The curve described by $\la_1$ is at height $\sqrt{2}$ exactly for $\tht=\pi/2+\om$, and for $\tht$
	in $[\pi/2-\om,\pi/2+\om]$, the point $\la_0=\sqrt{2}$ lies in $\si(A+B)$ and gives $\vro([A,B])=
	\sqrt{2(4-2)}=2$. For $\pi/2+\om <\tht < \pi-\om$, we have values $c-\ga$ whose square is closer to $2$
	than $(-\la_1)^2$; in the left picture of Figure~4, these values occur for $\tht$ from about $2.0$
	to about $2.8$. However, for $\tht \in (\om,\pi-\om)$  we have $s > \de$, and hence, by Theorem~9, 
	these points $c-\ga$ do not belong to $\si(A+B)$. Thus, for $\pi/2+\om < \tht \le \pi-\om$, the points
	$\la_1$ and $-\la_1$ are those in $\La$ whose squares are closest to $2$. It follows that $\vro([A,B])= 
	\sqrt{\la_1^2(4\!-\!\la_1^2)}$ for these $\tht$, which, as above, gives $\vro([A,B])=2|\sin(\tht-\om)|$.
	Let finally $\tht$ be in $(\pi-\om,\pi)$. In this case $c-\ga > -\la_1$, implying that the square of
	$c-\ga$ is closer to $2$ than the square of $-\la_1$, but from Lemma~8 we know that $c-\ga$ 
	is never in $\si(A+B)$, so that
	again $\vro([A,B])=\sqrt{\la_1^2(4\!-\!\la_1^2)}$.
	
	\smallskip
	The case where $\pi/2 \le \om < \pi$ can be disposed of similarly. \;\: $\square$

	% ------------------------------------------------------------------------
	
	%\subsection*{Acknowledgment}
	%Many thanks to our \TeX-pert for developing this class file.

	% ------------------------------------------------------------------------

\begin{thebibliography}{9}
		
		\bibitem{Uni}
		A. B\"ottcher and B. Silbermann:
		{\em Introduction to Large Truncated Toeplitz Matrices.}
		Universitext, Springer-Verlag, New York 1999.
		
		\bibitem{BS}
		A. B\"ottcher and I. M. Spitkovsky:
		{\em A gentle guide to the basics of two projections theory.}
		Linear Algebra Appl. 432, 1412--1459 (2010).
		
		\bibitem{RS}
		A. B\"ottcher and I. M. Spitkovsky:
		{\em Robert Sheckley's Answerer for two orthogonal projections.}
		Operator Theory: Adv. and Appl. 268, 125--138 (2018).
		
		\bibitem{CB}
		A. Chefles and S. M. Barnett:
		{\em  Diagonalisation of the Bell-CHSH operator. }
		Physics Letters A, vol. 232, no. 1-2, 4--8 (1997).
		
		\bibitem{CHSH}
		J. F. Clauser, M. A. Horne, A. Shimony, and R. A. Holt:
		{\em Proposed experiment to test local hidden-variable theories.}
		Physical Review Letters 23, 880--884 (1970).
		
		\bibitem{FT}
		Y. Fujii and T. Tsurumaru:
		{\em Estimating spectral radius of Bell-type operator via finite dimensional 
			approximation of orthogonal projections.}
		arXiv:2511.10939 [math.RT] (2025).
		
		\bibitem{GRS}
		M. S. Guimaraes, I. Roditi, and S. P. Sorella:
		{\em Introduction to Bell's inequality in Quantum Mechanics.}
		arXiv:2409.07597v1 [quant-ph] (2024).
		
		\bibitem{H}
		P. R. Halmos: 
		{\em Two subspaces.} 
		Trans. Amer. Math. Soc. 144, 381--389 (1969).
		
		
		\bibitem{KT} 
		L. A. Khalfin and B. S. Tsirelson: 
		{\em Quantum and quasi-classical analogs of
			Bell inequalities.} In: Symposium on the Foundations of Modern Physics,
		World Scientific Publishing, 441--460 (1985).
		
		\bibitem{L}
		L. J. Landau:
		{\em  Experimental tests of general quantum theories.}
		Letters in Mathematical Physics 14, 33--40 (1987).
		
		\bibitem{LS}
		G. S. Litvinchuk and I.M. Spitkovsky:
		{\em Factorization of Measurable Matrix Functions.}
		Birkh\"auser Verlag,  Basel and Boston 1987. 
		
		
		\bibitem{T}
		B. S. Tsirelson:
		{\em Quantum generalizations of Bell's inequality.}
		Letters in Mathematical Physics 4, 93--100 (1980).
		
	\end{thebibliography}
\end{document}